\documentclass[12pt,leqno]{article}
\usepackage[francais]{babel}
\usepackage[latin1]{inputenc}
\usepackage[T1]{fontenc}
\usepackage{amsmath}
\usepackage{amsfonts}
\usepackage{amssymb}

\newcommand{\Z}{\mathbb{Z}}
\renewcommand{\Re}{\mathop{\rm Re}}

\newcommand{\Ff}{\mathcal{F}}

\newcommand{\ba}{\backslash}
\newcommand{\rg}{\rightarrow}

\begin{document}

\title
{
Analyse Harmonique\\
Idempotents et échantillonnage parcimonieux\\
Idempotents and compressive sampling}

\author{Jean--Pierre Kahane}

\date{}

\maketitle

\noindent \textbf{Résumé}.

Comment reconstituer un signal, assimilé à une fonction $x$ définie sur le groupe cyclique $\Z_N$, qu'on sait porté par $T$ points, en n'utilisant sa transformée de Fourier $\hat{x}$ que sur un ensemble $\Omega$ de fréquences ? Le procédé indiqué par Candès, Romberg et Tao \cite{Cand,CaRoTao} est l'extrapolation minimale de $\hat{x}|_\Omega$ dans~$\Ff \ell^1$. La note traite les questions suivantes :
1) Quand est--il vrai que ce procédé redonne tous les signaux portés par $T$ points ? 2) Si l'on choisit $\Omega$ par sélection aléatoire de points de $\Z_N$, $N$ étant très grand, avec quelle probabilité obtient--on par ce procédé tous les signaux portés par $T$ points ? 3) tous les signaux portés par un ensemble $S$ donné ? 4) un signal donné ? Je donne des réponses à 1) et à 2) avec démonstrations, et à 3) sans démonstration. La réponse à 3) améliore les estimations de Candès, Romberg et Tao relatives à 4), la question qu'ils traitent. L'idempotent $K$ tel que $\hat{K}=1_\Omega$ joue un rôle central.

\vskip2mm

\noindent \textbf{Abstract}.

According to Candès, Romberg and Tao \cite{Cand,CaRoTao}, a signal is represented as a function $x$ defined on the cyclic group $\Z_N$. Assuming that it is carried by a set $S$ consisting of $T$ points, how to reconstruct $x$ by using only a small set $\Omega$ of frequences ? The procedure of Candès, Romberg and Tao is the minimal extrapolation of $\hat{x}|_\Omega$ in $\Ff\ell^1$, when it exists. 1) When can we obtain in this way all signals carried by $T$ points ? 2) Choosing $\Omega$ by a random selection of points in $\Z_N$ with $N$ very large, give an estimate of the probability that the procedure works for all signals carried by $T$ points 3) for all signals carried by a given set $S$ 4) for a given signal. The answers to 1) and 2) are given with proofs and the answer to 3) without proof. Candès, Romberg and Tao answered question 4) and our answer to 3) improves their estimates. A key role is played by the idempotent $K$ such that $\hat{K}=1_\Omega$.

\vskip 4mm

Cette note est inspirée d'un théorème de Candès, Romberg et Tao (Theorem~1.3 de \cite{CaRoTao}, Theorem~2.1 de \cite{Cand}). Ce théorème donne une méthode pour reconstruire un signal porté par $T$ points de $\Z_N$ à partir de la restriction de sa transformée de Fourier à un ensemble $\Omega$ de fréquence bien plus petit que $\Z_N$. C'est le modèle typique d'un échantillonnage parcimonieux. Voici le cadre et l'énoncé.

$\Z_N$ est le groupe $G$ des temps $t$, et aussi bien le groupe $\hat{G}$ des fréquences~$\omega$ ; la dualité s'exprime~par
$$
\langle\omega,t\rangle = e\Big(\frac{\omega t}{N}\Big)\,, \ e(\xi) = e^{2 \pi i\xi}\,.
$$
Le signal est une fonction $x(t)$ $(t\in G)$ à valeurs complexes, sa transformée de Fourier est
$$
\hat{x}(\omega) = \frac{1}{\sqrt{N}} \sum_{t\in G} x(t) e \Big(\frac{-\omega t}{N}\Big)
$$
et la reconstruction classique se fait par la formule d'inversion
$$
x(t) = \frac{1}{\sqrt{N}} \sum_{\omega\in \hat{G}} \hat{x}(\omega) e\Big(\frac{\omega t}{N}\Big)\,.
$$
Soit $A(\hat{G})$ l'image de Fourier de $\ell^1 (G)$ :
$$
\| \hat{x}\|_{A(\hat{G})} = \|x\|_{\ell^1(G)} = \sum_{t\in G} |x(t)|\,.
$$
Sous certaines conditions, qu'on va expliciter,

\vskip2mm

$(i)$ $\hat{x}$ \textit{est le prolongement minimal de $\hat{x}|_\Omega$ dans $A(\hat{G})$.}

\vskip2mm

La reconstruction se ramène alors à un problème d'extremum en analyse convexe, qui est praticable.

\vskip2mm

\noindent \textsc{Théorème CRT}.~--- \textit{Soit $x$ un signal porté par une partie $S$ de $G$, de cardinal $T:|S|=T$. Soit $\Omega$ une partie aléatoire de $\hat{G}$, dont la distribution de probabilité est uniforme sur l'ensemble des parties de $\hat{G}$ de cardinal $|\Omega|$ donné~par}
$$
|\Omega| = [CT\log N]\quad ([\ ] : \textit{partie\ entière})
$$
\textit{avec}
$$
C = 22 (1+\delta)\,, \quad \delta >0\,.
$$
\textit{Alors la probabilité de l'évènement (i) vérifie}
$$
P((i)) = 1 - O(N^{-\delta}) \quad (N \rg \infty) \qquad \cite{Cand}
$$
\vskip2mm

Pour des $N$ et $\delta$ convenables, $(i)$ est très probable. Quitte à changer la définition de $\Omega$, peut--on dire qu'il est très probable que $(i)$ ait lieu pour tous les $x$ portés par un $S$ donné ? ou mieux, pour tous les $x$ portés par un $S$ de cardinal $T$ donné ? ou mieux encore, peut--on dire que cela est certain ? Je vais proposer quelques réponses.

D'abord, comme il est noté dans \cite{Cand}, $(i)$ est entraîné par

\vskip2mm

 $(ii)$ \textit{pour tout $z\in \ell^1(G)$ non nul, tel que $ \hat{z}|_\Omega=0$,}
$$
\sum_{t\in G\ba S}|z(t)| > \sum_{t\in S} |z(t)|\,. 
$$
On voit facilement que $(ii)$ est entraîné par

\vskip2mm

$(iii)$ \textit{pour toute fonction $\lambda$ de module $1$ sur $S$, il existe un $p \in \ell(G)$\break tel que $\hat{p}$ soit porté par $\Omega$ et~que}
$$
\left\{
\begin{array}{lll}
|p(t)-\lambda(t)| < \dfrac{1}{2} &\mathrm{quand} &t\in S\\
\noalign{\vskip2mm}
|p(t)| < \dfrac{1}{2} &\mathrm{quand} &t\in G\ba S
\end{array}
\right.
$$

Faisons intervenir l'idempotent
$$K(t) = \sum_{\omega \in \Omega} e\Big(\frac{\omega t}{N}\Big)\,,
$$
et cherchons $p$ sous la forme
$$
p(t) = \sum_{t' \in S} \lambda (t') \frac{K(t-t')}{K(0)}\,.
$$
On voit que $(iii)$ est vérifié sous la condition que

\vskip2mm

$(iv)$ \textit{pour tout $t\not= 0$},
$$
|K(t)| < \frac{1}{2T} K(0)\,.
$$

Voici donc un résultat certain.

\vskip2mm

\textbf{T1.} \textit{Si l'idempotent $K(t) = \sum\limits_{\omega\in \Omega}\langle\omega,t\rangle$ vérifie $(iv)$, $(i)$ a lieu pour tous les signaux $x$ portés par un ensemble de $T$ points.} 

Ce qui importe dans cet énoncé est que le spectre de $K$ soit contenu dans~$\Omega$.  Les idempotents apparaissent de façon naturelle dans la construction qui
suit. En effet, reste à construire $\Omega$ de façon que $(iv)$ ait lieu. Pour cela, on reprend la construction de \cite{CaRoTao}. Au lieu de fixer $|\Omega|$, on fixe un $0<\tau<1$, on considère des variables aléatoires de Bernoulli indépendantes d'espérance $\tau$ \hbox{$(P(X_n=1)=\tau$}, $P(X_n =0)=1-\tau)$ $(n\in \Z_N)$ et on pose $\Omega = \{ n|X_n =1\}$. Ainsi
$$
K(t) = \sum_{n \in \Z_N} X_n e \Big(\frac{\omega t}{N}\Big)\,.
$$
Dans les énoncés qui suivent, $\tau N$ jouera le rôle de $|\Omega|$.

Majorons la probabilité pour que $(iv)$ n'ait pas lieu. Fixons $t\not= 0$ et un entier $\nu>3$. Si $|K(t)| \ge \dfrac{1}{2T} K(0)$,  il existe un $\varphi= \frac{2j\pi}{\nu}$ tel que
$$
\Re K(t) e^{-i\varphi} \ge \cos \frac{\pi}{\nu} \,\frac{1}{2T}K(0)\,.
$$
Posons $a=\cos \frac{\pi}{\nu}$ et
$$
\begin{array}{ll}
Y &=\Re K(t) e^{-i\varphi} - \dfrac{a}{2T} K(0)\\
\noalign{\vskip2mm}
&= \displaystyle\sum_{n\in \Z_N} X_n \Big(\cos \Big(\frac{2\pi nt}{N} -\varphi\Big) - \frac{a}{2T}\Big) = \sum_{n\in \Z_N}X_n A_n\,.
\end{array}
$$
Pour $u>0$,
$$
\begin{array}{c}
P(Y\ge 0) < E e^{uY} = \displaystyle\prod_{n\in\Z_N} E e^{uX_nA_n}\\
\noalign{\vskip2mm}
= \displaystyle \prod_{n\in \Z_N} (1-\tau +\tau e^{uA_n}) \le \prod_{n\in \Z_N} \exp (\tau(e^{uA_n}-1)\,.
\end{array}
$$
Or
$$
\sum_{n\in Z_N} (e^{uA_n}-1) = \exp \Big(-\frac{au}{2T}\Big) \sum_{n\in \Z_N} \exp \Big(u \cos \Big(\frac{2\pi nt}{N}-\varphi\Big)\Big) -N\,.
$$
En développant en puissances de $u$ on est amené à choisir $u=\frac{a}{T}$, d'où
$$
\begin{array}{lcl}
P(Y\ge 0) &<& \exp \Big(\tau N \Big(-\dfrac{a^2}{4T^2} + \dfrac{a^4}{8T^4} + \displaystyle \sum_{k=4} \dfrac{a^k}{k!T^k}\Big)\\
\noalign{\vskip2mm}
&<& \exp \Big(TN \Big( -\dfrac{a^2}{4T^2} + \dfrac{a^4}{4T^4}\Big)\Big)
\end{array}
$$
sous les conditions $t\not=0$ (hypothèse déjà faite), $2t\not=0$ et $3t\not=0$. Donc
$$
1-P ((iv)) < N\nu \exp \Big(\tau N \Big( - \frac{a^2}{4T^2} + \frac{a^4}{4T^4}\Big)\Big) \qquad \Big(a=\cos \frac{\pi}{\nu}\Big)\,. \leqno(1)
$$

Choisissons
$$
\tau N = 4CT^2 \log N\,. \leqno(2)
$$

Le second membre de (1) devient $\nu N^{-d}$ avec
$$
d=C \Big(a^2 - \frac{a^4}{T^2}\Big) -1\,.
\leqno(3)
$$

\textbf{T2.} \textit{Supposons $N=\pm 1$ modulo~$6$. La condition $(iv)$, donc la condition (i) pour tous les signaux $x$ portés par $T$ points de $\Z_N$, est réalisée avec une probabilité supérieure à $1-\nu N^{-d}$ quand on choisit $\Omega$ selon $(2)$, dès que $\nu,C$ et $d$ sont liés par $(3)$ avec $a=\cos \frac{\pi}{\nu}$.}

Par exemple, si $T=2$, $C=2$, $N=1001$ et $\nu=10$, on a une probabilité supérieure à $\frac{1}{3}$ de réaliser $(iv)$, donc $(i)$ pour tous les signaux portés par deux points, en prenant au hasard $\Omega$ dans $\Z_N$ avec $|\Omega|=222$. En soignant les calculs, on peut remplacer $1/3$ par~$3/4$. En tout cas, ce résultat n'est pas glorieux.

Il est bon de remarquer que le facteur $4T^2$ dans (2) est dans la nature de la question : si $(iv)$ est réalisé pour $K(t)=\sum\limits_{\omega\in \Omega}\langle\omega,t\rangle$, il est facile de voir~que
$$
|\Omega| \ge 4T^2 \frac{N}{N+4T^2-1}\,. \leqno(4)
$$
Par exemple, si $T=2$ et $N=1001$, on a $|\Omega| \ge 16$.

\vskip2mm

\textit{Voici un corollaire de $T2$ : dès que}
$$
C > \frac{T^2}{T^2-1}\,,
\leqno(5)
$$
\textit{le choix de $\Omega$ selon $(2)$ entraîne que, pour $N=\pm 1$ modulo $6$ assez grand, la procédure décrite par $(i)$ reconstruit tous les signaux portés par $T$ points avec une probabilité aussi voisine de $1$ que l'on veut.}

\vskip2mm

En travaillant un peu plus, on peut retrouver en l'améliorant le théorème CRT sous la forme suivante, qui répond à la première question que nous avons posée :

\vskip2mm

\textbf{T3.} \textit{La condition $(ii)$, donc la condition $(i)$ pour tous les signaux $x$ portés par un ensemble $S \subset \Z_N$, de cardinal $T$, est réalisée avec une probabilité supérieure à $1-h$ quand on choisit $\Omega$ selon}
$$
\tau N = 4CT \log N
\leqno(6)
$$
\textit{et qu'on prend}
$$
h \ge T \mu^T \nu N^{-Ca^2\alpha'^2} + (N-T) \mu^T \nu N^{-Ca^2(1-\alpha)^2} 
\leqno(7)
$$
\textit{$\mu$, $\nu$, $a$, $\alpha$ et $\alpha'$ étant soumis aux conditions suivantes} : $\mu$ \textit{et $\nu$ entiers} $>1$, $a=\cos \frac{\pi}{\nu}$, $0<\alpha<1$, $\alpha' =\alpha -2 \sin \frac{\pi}{2\mu}$.

Par exemple, en prenant $\mu=10$, $\nu=10^3$, $\alpha=\frac{2}{3}$, on peut choisir
$$
h = 10^{T+3} N^{-\frac{C}{10}+1}\,,
\leqno(8)
$$
ce qui est explicite mais sans aucun intérêt pratique. Si on désire une estimation en $O(N^{-\delta})$ $(N\rg \infty)$ comme dans le théorème~CRT, on choisit $\mu$ et $\nu$ très grands et $\alpha$ convenable, et on obtient :

\vskip2mm

\noindent \textsc{Corollaire de T3}.~--- \textit{Dès que $C>1$, le choix de $\Omega$ selon $(6)$ entraîne au lieu de $(7)$ l'évaluation}
$$
h= O(N^{-\delta}) \ (N\rg \infty)
\leqno(9)
$$
\textit{dès que}
$$
\delta < \frac{(C-1)^2}{4C}\,.
\leqno(10)
$$

Cela se compare avantageusement à l'estimation de CRT.

La démonstration de T3 sera donnée par ailleurs, ainsi que d'autres variations possibles sur le théorème CRT.

\eject

\end{document}